\documentclass[twoside,12pt]{article}
\usepackage{amsfonts,amsmath,amsthm,amssymb,amscd}
\title {On the Jacobi Group and the\\
Mapping Class Group of $S^3\times S^3$}
\author{Nikolai A. Krylov
\\
\small University of Illinois at Chicago \\
\small Department of Mathematics, Statistics and Computer Science\\
\small 851 S.Morgan st. Chicago, IL 60607\\
\small e-mail:~krylov@math.uic.edu}

\begin{document}

\date {}

\newtheorem{thm}{Theorem}
\newtheorem{lem}{Lemma}
\newtheorem{claim}{Claim}
\newtheorem{dfn}{Definition}

\def\defin          {\stackrel{\it def}{=}}
\def\mcgss          {\pi_0 {\it Diff(S^3\times S^3)} }
\def\mcgM           {\pi_0 {\it Diff(M)} }
\def\kerM           {\pi_0S {\it Diff(M)} }
\def\kerss          {\pi_0S {\it Diff(S^3\times S^3)} }
\def\zint           {\mathbb Z}
\def\faczin         {\mathbb Z_{28}}
\def\lra            {\longrightarrow}
\def\ra             {\rightarrow}
\def\hra            {\hookrightarrow}
\def\lmt            {\longmapsto}
\def\G1             {SL_2(\mathbb Z)}
\def\Aut            {Aut~H_k(M)}
\def\vars           {\varsigma}
\def\comAB          {ABA^{-1}B^{-1}}
\def\cH             {{\cal H}}
\def\cM             {{\cal M}}
\def\ZZ             {\zint\oplus\zint}
\def\GJ             {\Gamma^J}
\def\lam            {\lambda}
\def\del            {\delta}
\def\to             {\tilde\omega}
\def\wS             {\widetilde{\Sigma}}
\def\wH             {\widehat{\Sigma}}

\maketitle


\begin{abstract}
{The paper contains a proof that the mapping class group of the
manifold $S^3\times S^3$ is isomorphic to a central extension of
the (full) Jacobi group $\GJ$ by the group of 7-dimensional
homotopy spheres. Using a presentation of the group $\GJ$ and the
$\mu$-invariant of the homotopy spheres, we give a presentation
of this mapping class group with generators and defining
relations. We also compute cohomology of the group $\GJ$ and
determine a 2-cocycle that corresponds to the mapping class group
of $S^3\times S^3$.}
\end{abstract}

\section{Introduction}

The central theme of this paper is the group of isotopy classes
of orientation preserving diffeomorphisms on $S^3\times S^3$. We
will denote this group by $\mcgss$. In general, the group of
isotopy classes of orientation preserving diffeomorphisms on a
closed oriented smooth manifold $M$ will be denoted by $\mcgM$
and called the {\it mapping class group} of $M$ by analogy with
the 2-dimensional case.

The article consists of two parts. Our goal in the first part
will be to give a presentation of the mapping class group of
$S^3\times S^3$ with generators and defining relations. The main
step in this direction is Theorem 1. where we prove that $\mcgss$
is a central extension of the (full) Jacobi group $\GJ$ by the
group of 7-dimensional homotopy spheres $\Theta_7$. The second
part is concerned with the cohomology group
$H^2(\GJ,\zint_{28})$. We show that this group is isomorphic to
$\zint_{28}\oplus\zint_4\oplus\zint_2$ and determine a 2-cocycle
that corresponds to $\mcgss$.

\noindent {\bf Acknowledgements}\\
Results of this paper are part of the author's doctoral thesis
written under the direction of Anatoly Libgober. I am deeply
indebted to him for advice and continuing encouragement.

\section{A presentation of the group $\mcgss$}

For a closed smooth oriented (k-1)-connected
almost-parallelizable manifold $M^{2k},~k\geq2$ the group $\mcgM$
has been computed in terms of exact sequences by Kreck
\cite{Kreck}. We will begin by recalling his results in the first
section. Next we will show that there is an isomorphism between
the factor group $\left.\mcgss\right/\Theta_7$ and the full
Jacobi group $\GJ$. The group of isotopy classes of
diffeomorphisms that act trivially on homology of $S^3\times S^3$
is determined in $\S 2.3$. In the last section of this part we
give a presentation of $\mcgss$ with generators and defining
relations (Theorem 3.). All diffeomorphisms are assumed to be
orientation preserving and integer coefficients are understood
for all homology and cohomology groups.

\subsection{Exact sequences of Kreck}

Our focus in this paragraph is to review the results of Kreck
\cite{Kreck}. First we recall some of the definitions and
notations. From now on by a manifold $M$ we will mean a closed
oriented differentiable (k-1)-connected almost-parallelizable
2k-manifold.

Denote by $\Aut$ the group of automorphisms of $H_k(M,\zint)$
preserving the intersection form on $M$ and (for $k\geq3$)
commuting with the function $\alpha:~H_k(M)\lra
\pi_{k-1}(SO(k))$, which is defined as follows. Represent $x\in
H_k(M)$ by an embedded sphere $S^k\hookrightarrow M$. Then
function $\alpha$ assigns to $x$ the classifying map of the
corresponding normal bundle. Any diffeomorphism $f\in {\it
Diff(M)}$ induces map $f_*$ which lies in $\Aut$. This gives a
homomorphism
$$
\kappa:~ \mcgM \lra \Aut,~~~[f]\longmapsto f_*
$$

The kernel of $\kappa$ is denoted by $\kerM$. For elements of
$\kerM$ Kreck defines the following invariant:~ Choose again a
sphere $S^k\hookrightarrow M$ that represents an element  $x\in
H_k(M)$. Since $[f]\in \kerM$ we can assume that $f|_{S^k}=Id$.
The stable normal bundle $\nu(S^k)\oplus1$ of $S^k$ in $M$ is
trivial and we can choose some trivialization
$\tau:~\nu(S^k)\oplus 1 \lra S^k\times \mathbb R^{k+1}$. Clearly,
the differential of $f$ leaves the tangent bundle of $S^k$
invariant and hence induces an automorphism of the normal bundle
$\nu(S^k)\oplus1$. At each point $t\in S^k$ this automorphism
gives (via trivialization $\tau$) an element ${\cal P}_t$ of the
group $SO(k+1)$ and hence we get an element ${\cal P}\in
\pi_k(SO(k+1))$. It is obvious that $\cal P$ lies in the image of
the map~ $S:~\pi_k(SO(k))\lra\pi_k(SO(k+1))$~ induced by the
inclusion $SO(k)\hookrightarrow SO(k+1)$. It is a standard fact
that element $\cal P$ does not depend on trivialization $\tau$ of
the normal bundle $\nu(S^k)\oplus1$. This construction leads to a
well defined homomorphism (cf. Lemmas 1,2 of \cite{Kreck})
$$
\chi:~\kerM\lra Hom(H_k(M),~S\pi_k(SO(k)))$$

For $k\equiv 3\pmod{4}$ the group $S\pi_k(SO(k))$ is isomorphic
to the cyclic group $\zint$ and hence we can identify
$Hom(H_k(M),~S\pi_k(SO(k)))$ with the cohomology group $H^k(M)$.
In this case one can describe $\chi(f)$ by the Pontrjagin class
of the mapping torus $M_f$, but we will not use this description
here so we omit the details. The following theorem is due to
Kreck (cf. {\underline{Theorem 2}}, \cite{Kreck}):\\

\noindent {\bf Theorem:} $k\geq 3$. If $M^{2k}$ bounds a framed
manifold, then the following sequences are exact:
\begin{equation}
\label{firstseq} 0\lra \kerM \lra \mcgM \stackrel{\kappa}{\lra}
\Aut \lra 0
\end{equation}
\begin{equation}
\label{secondseq} 0\lra \Theta_{2k+1} \stackrel{\iota}{\lra} \kerM
\stackrel{\chi}{\lra} Hom(H_k(M),~S\pi_k(SO(k))) \lra 0
\end{equation}\\
\noindent {\underline{Remark:}} If one considers $M^{2k}$ which
does not bound a framed manifold then ${\it Ker}(\chi)$ will be a
factor group $\left.\Theta_{2k+1}\right/\Sigma_M$ instead of the
whole group of homotopy spheres. Since $S^3\times S^3$ is the
boundary of $S^3\times D^4$, we have $\Sigma_M=0$ (cf. Lemma 3.
of \cite{Kreck}).

Map $\iota$ is defined as follows. Present a homotopy sphere
$\Sigma\in\Theta_{2k+1}$ as union $D^{2k+1}\bigcup_f D^{2k+1}$ and
assume that $f=id$ on a neighbourhood of the lower hemisphere
$D^{2k}_-\subset S^{2k}$. Then $\iota(\Sigma)$ is the class of
diffeomorphism on $M$ which is identity outside an embedded disk
in $M$ and is equal to $f|_{D^{2k}_+}$ on this disk.

It follows from the second exact sequence that
$\left.\kerss\right/\Theta_7$ is isomorphic to $\ZZ$. Note also
that $Aut~H_3(S^3\times S^3)\cong\G1 $. Indeed, the function
$\alpha: H_3(S^3\times S^3)\lra \pi_2(SO(3))$ is zero map. Hence
an element $A$ of the group $Aut~H_3(S^3\times S^3)$ will be any
automorphism of $\ZZ$ which preserves the intersection form on
$S^3\times S^3$. It means that if we choose a basis for
$H_3(S^3\times S^3)$ then $A\in\G1 $. Exact sequences
(\ref{firstseq}) and (\ref{secondseq}) induce the following short
exact sequence
\begin{equation}
\label{thirdseq} 0\lra \ZZ \lra \left.\mcgss\right/\Theta_7 \lra
\G1  \lra 0
\end{equation}

In the next section we prove that this exact sequence splits and
the group $\left.\mcgss\right/\Theta_7$ is isomorphic to the
Jacobi group $\GJ$.

\subsection{Splitting of the Exact Sequence}

The (full) Jacobi group $\GJ$ is a semidirect product of the
modular group with the direct sum $\ZZ$. More precisely
(cf. \cite{Zagier} \S I.1.),\\
$\Gamma^J\defin \G1 \ltimes\zint^2$ =~set of pairs $(M,X)$ with
$M\in \G1 ,~X\in\zint\oplus\zint$ and group law
$(M,X)\cdot(M',X')=(MM',XM'+X')$ (notice that vectors are written
as row vectors, i.e. $\G1 $ acts on the right). It is interesting
to note that $\GJ$ first came up in the theory of
Jacobi forms (see \cite{Zagier}). We will need
\begin{lem}
$\Gamma^J$ admits the following presentation:\\
$<y,u,a,b~|~yuy=uyu,~(yuy)^4=id,~ab=ba,~ay=yab,~au=ua,~
by=yb,$\\ $bu=uba^{-1}>$
\end{lem}
\begin{proof}
$\G1 $ has a presentation: $<y,u~|~yuy=uyu,~(yuy)^4=id>$ (see for
example \cite{Bir.paper}) where $y$ and $u$ correspond to
matrices $\begin{pmatrix} 1&1\\0&1\end{pmatrix}$,
$\begin{pmatrix} 1&0\\-1&1\end{pmatrix}$ respectively. ({\it It is
a classical fact that $\G1 \cong\zint_4*_{\zint_2}\zint_6$. Hence
$\G1 $ has a presentation: $<x,z~|~x^4=id,~x^2=z^3>$. One can use
a map $f$: $f(z)=yu,~f(x)=(yuy)^{-1}$ to show that these two
presentations define isomorphic groups.}) By definition of
$\Gamma^J$ the following sequence is exact.
$$ 0\lra \zint^2 \lra \Gamma^J \lra
\G1 \lra 0$$

Consider a homomorphism $\alpha: \G1 \lra \Gamma^J$ defined by
the formulas:
$\alpha(y)\defin(y,(0,0)),~~\alpha(u)\defin(u,(0,0))$. If we
denote elements $(y,(0,0))$; $(u,(0,0))$;
$(id,(1,0));~(id,(0,1))$ by $y,u,a$ and $b$ respectively we see
that these elements $y,u,a,b$ generate $\Gamma^J$ and the
relations $yuy=uyu,~(yuy)^4=id,~ab=ba$ are satisfied. To find all
defining relations for $\Gamma^J$ we need to find how $\G1 $ acts
on the generators $a$ and $b$ of $\zint^2$ by conjugation. First
note that $ab=(id,(1,0))\cdot(id,(0,1))=(id,(1,1))$ and\\
$ba^{-1}=(id,(0,1))\cdot(id,(-1,0))=(id,(-1,1))$. Hence\\
$ay=(id,(1,0))\cdot(y,(0,0))=(y,(1,0)\cdot y+(0,0))=(y,(1,1))$\\
$yab=(y,(0,0))\cdot(id,(1,1))=(y,(1,1))\Rightarrow ay=yab$\\
$au=(id,(1,0))\cdot(u,(0,0))=(u,(1,0)\cdot u+(0,0))=(u,(1,0))$\\
$ua=(u,(0,0))\cdot(id,(1,0))=(u,(1,0))\Rightarrow au=ua$\\
$by=(id,(0,1))\cdot(y,(0,0))=(y,(0,1)\cdot y+(0,0))=(y,(0,1))$\\
$yb=(y,(0,0))\cdot(id,(0,1))=(y,(0,1))\Rightarrow by=yb$\\
$bu=(id,(0,1))\cdot(u,(0,0))=(u,(0,1)\cdot u+(0,0))=(u,(-1,1))$\\
$uba^{-1}=(u,(0,0))\cdot(id,(-1,1))=(u,(-1,1))\Rightarrow
bu=uba^{-1}$
\end{proof}

\noindent {\underline{Remark:}} A different presentation of this
group can be found in \cite{Choie} (cf. also Lemma 6. below).

Consider now the standard sphere $S^3$ in Euclidean four-space
$\mathbb R^4$, given by the equation:
$x_0^2+x_1^2+x_2^2+x_3^2=1$. This sphere can be identified with
the special unitary group $SU(2)$ which is also known as group of
unit quaternions. The group structure on $S^3$ induces the group
structure on the product $S^3\times S^3$. If we denote elements of
the group $S^3$ by $s$ and $t$ we will write $(s,t)$ to denote the
corresponding element of the group $S^3\times S^3$. The product
of two elements $(s,t)$ and $(s',t')$ will be the pair
$(ss',tt')$ with quaternion multiplication understood.

\begin{thm}
The factor group $\pi_0 {\it Diff(S^3\times S^3)}/{\Theta_7}$ is
isomorphic to the Jacobi group $\Gamma^J$.
\end{thm}

\begin{proof}
We will give a presentation of the factor group $\pi_0 {\it
Diff(S^3\times S^3)}/{\Theta_7}$ that coincides with the above
presentation of $\Gamma^J$. By $y$ and $u$ we denote
the generators of $\G1 $ as above. Consider isotopy classes $[Y]$ and
$[U]$ of the following diffeomorphisms of $S^3\times S^3~
((s,t)\in S^3\times S^3)$:
\begin{equation}
\label{twists} Y: (s,t)\longmapsto (s,st),~~~~U: (s,t)\longmapsto
(t^{-1}s,t)
\end{equation}
Define map $\beta: \G1 \ra \mcgss$ by the identities:
$\beta(y)\defin [Y]$, $\beta(u)\defin [U]$ and extend it linearly
to the whole group $\G1 $. We will show that $\beta$ is a well
defined homomorphism from $\G1 $ to $\pi_0 {\it Diff(S^3\times
S^3)}/{\Theta_7}$.
First we check that $YUY=UYU$:\\
$YUY:~~~(s,t)\stackrel{Y}{\longmapsto}(s,st)\stackrel{U}{\longmapsto}
((st)^{-1}s,st)=(t^{-1},st)\stackrel{Y}{\longmapsto}(t^{-1},t^{-1}st)$\\
$UYU:~~~(s,t)\stackrel{U}{\longmapsto}(t^{-1}s,t)\stackrel{Y}{\longmapsto}
(t^{-1}s,t^{-1}st)\stackrel{U}{\longmapsto} (t^{-1},t^{-1}st)$\\
Thus $YUY=UYU$ and hence $[Y][U][Y]=[U][Y][U]$. From now on we
will denote a diffeomorphism and the isotopy class of it by the
same capital letter, omitting the brackets. To prove the equality
$(YUY)^4=Id$~ ($Id$ stands for the identity diffeomorphism of
$S^3\times S^3$)
we will need some auxiliary results. Consider the following
diffeomorphisms $A,B\in {\it Diff(S^3\times S^3)}$:
\begin{equation}
\label{defAB}
A:~(s,t)\longmapsto(tst^{-1},t)~~~~~~
B:~(s,t)\longmapsto(s,sts^{-1})
\end{equation}
If we choose spheres $S^3\times1$ and $1\times S^3$ as generators
of the group $H_3(S^3\times S^3)$, it is obvious that
diffeomorphisms $A$ and $B$ preserve these spheres and act
trivially on homology of $S^3\times S^3$.
\begin{lem}
Isotopy classes of diffeomorphisms $A$ and $B$ generate the
group\\ ${\it Hom(H_3(S^3\times S^3),~S\pi_3(SO(3)))}\cong
\zint\oplus\zint$
\end{lem}
\begin{proof}[Proof of the lemma]
Let us compute $\chi(B)$. Since $S^3\times S^3$ is a
parallelizable manifold and the normal bundle of $S^3_1\defin
S^3\times1$ in $S^3\times S^3$ is trivial we need to find an
element of the group $\pi_3(SO(3))$ that corresponds to the
differential of $B$. Take a point $(s,1)\in S^3_1$. We can
identify the fiber of the normal bundle $\nu(S^3_1)$ at $(s,1)$
with the fiber of the tangent bundle $\tau(s\times S^3)$ at this
point. Furthermore, via the projection $\rho_2: S^3\times S^3\lra
1\times S^3$ we can identify this tangent fiber at $(s,1)$ with
the tangent fiber at $(1,1)$ (Lie algebra $\mathfrak g$ of the
group $1\times S^3$). Since $B_s(t)\defin B|_{s\times S^3}(s,t)=
sts^{-1}$ the map $s\mapsto d_{(s,1)}B_s$ will correspond to the
{\it adjoint representation} $Ad: S^3_1\lra Aut(\mathfrak g)$, and
$\chi(B)(s,1)=Ad_s\in SO(3)$. Thus
$\chi(B)(S^3\times1)=Ad:~S^3\times1\lra SO(3)$. If we choose an
element $T\in \mathfrak g$ then it is well known that
$Ad_s(T)=sTs^{-1}$. This map is a generator of the group
$\pi_3(SO(3))$ (see \cite{Lawson}, ch.I \S2) and therefore the
isotopy class of diffeomorphism $B$ is a generator of the group
$H^3(S^3\times S^3)$. In a similar way one can show that the
isotopy class of diffeomorphism $A$ is the other generator
(corresponding to the map $\chi(A)(1\times S^3):~1\times S^3\lra
SO(3)$) of the group $H^3(S^3\times S^3)$.
\end{proof}
Now we show that $AB=BA,~~AY=YAB,~~AU=UA$, $BY=YB,~~BU=UBA^{-1}$
in the factor group $\pi_0 {\it Diff(S^3\times S^3)}/{\Theta_7}$.
The first equality follows from the results of Kreck,
since as we just saw, $A$ and $B$ generate the abelian subgroup
of the group $\pi_0 {\it Diff(S^3\times S^3)}/{\Theta_7}$. To
prove the other equalities we use the group structure on
$S^3\times S^3$:\\
$AY:~~~~~(s,t)\stackrel{Y}{\longmapsto}(s,st)\stackrel{A}{\longmapsto}
(stst^{-1}s^{-1},st)$\\
$YAB:~~~(s,t)\stackrel{B}{\longmapsto}(s,sts^{-1})\stackrel{A}
{\longmapsto}(sts^{-1}sst^{-1}s^{-1},sts^{-1})=(stst^{-1}s^{-1},sts^{-1})$\\
$\stackrel{Y}{\longmapsto}
(stst^{-1}s^{-1},stst^{-1}s^{-1}sts^{-1})=(stst^{-1}s^{-1},st)~~
\Rightarrow AY=YAB,$\\
$AU:~~~~~(s,t)\stackrel{U}{\longmapsto}(t^{-1}s,t)\stackrel{A}
{\longmapsto}(tt^{-1}st^{-1},t)=(st^{-1},t)$\\
$UA:~~~~~(s,t)\stackrel{A}{\longmapsto}(tst^{-1},t)\stackrel{U}
{\longmapsto}(t^{-1}tst^{-1},t)=(st^{-1},t)~~\Rightarrow AU=UA,$\\
$BY:~~~~~(s,t)\stackrel{Y}{\longmapsto}(s,st)\stackrel{B}
{\longmapsto}(s,ssts^{-1})$\\
$YB:~~~~~(s,t)\stackrel{B}{\longmapsto}(s,sts^{-1})\stackrel{Y}
{\longmapsto}(s,ssts^{-1})~~\Rightarrow BY=YB,$\\
$BU:~~~~~(s,t)\stackrel{U}{\longmapsto}(t^{-1}s,t)\stackrel{B}
{\longmapsto}(t^{-1}s,t^{-1}sts^{-1}t)$\\
$UBA^{-1}:~(s,t)\stackrel{A^{-1}}{\longmapsto}(t^{-1}st,t)\stackrel{B}
{\longmapsto}(t^{-1}st,t^{-1}sttt^{-1}s^{-1}t)=(t^{-1}st,t^{-1}sts^{-1}t)$\\
$\stackrel{U}{\longmapsto}
(t^{-1}st^{-1}s^{-1}tt^{-1}st,t^{-1}sts^{-1}t)=(t^{-1}s,t^{-1}sts^{-1}t)~~
\Rightarrow BU=UBA^{-1}$\\ as required.
\begin{claim}
$(B^{-1}YUY)^4=Id,~~YUYB^{-1}=A^{-1}YUY,~~YUYA^{-1}=BYUY$.
\end{claim}
\begin{proof}[Proof of the claim:]
$B^{-1}YUY:~(s,t)\stackrel{YUY}{\lmt}(t^{-1},t^{-1}st)\stackrel{B^{-1}}
{\lmt}(t^{-1},tt^{-1}stt^{-1})=(t^{-1},s)$ and hence
$(B^{-1}YUY)^4:~(s,t)\lmt(t^{-1},s)\lmt(s^{-1},t^{-1})
\lmt(t,s^{-1})\lmt(s,t)$ i.e. $(B^{-1}YUY)^4=Id$.\\
Identities $AY=YAB,~BU=UBA^{-1},~AB=BA$ and $AU=UA$ imply
$YB^{-1}A^{-1}=A^{-1}Y,~~UB^{-1}=B^{-1}A^{-1}U$. Then from the
above equalities we get
$YUYB^{-1}=YUB^{-1}Y=YB^{-1}A^{-1}UY=A^{-1}YUY$. Similarly we see
that $YA^{-1}=BA^{-1}Y,~~UB=BUA$, hence
$YUYA^{-1}=YUBA^{-1}Y=YBUAA^{-1}Y=YBUY=BYUY$ which proves the
claim. \end{proof}

\noindent Using these identities we can show that $(YUY)^4=Id$.
$Id=(B^{-1}YUY)^4=$
$B^{-1}YUYB^{-1}YUYB^{-1}A^{-1}(YUY)^2=B^{-1}YUYB^{-1}BYUYB^{-1}(YUY)^2=$\\
$B^{-1}YUYYUYB^{-1}(YUY)^2=B^{-1}YUYA^{-1}(YUY)^3=(YUY)^4$. It
implies that exact sequence (\ref{thirdseq}) splits; the factor
group $\pi_0 {\it Diff(S^3\times S^3)}/{\Theta_7}$ has four
generators $Y,U,A,B$ and the following set of defining relations:
$YUY=UYU,~(YUY)^4=Id,~AB=BA,~AY=YAB,~AU=UA,~BY=YB,~BU=UBA^{-1}$.
In particular, we see that groups $\pi_0 {\it Diff(S^3\times
S^3)}/{\Theta_7}$ and $\Gamma^J$ have the same presentations and
therefore isomorphic.
\end{proof}
Note that diffeomorphisms $A$ and $B$ have been considered by
Browder (\cite{Browder}, Theorem 6.) to give an example of
diffeomorphisms of $S^3\times S^3$ which are homotopic to the
identity, but are not pseudo-isotopic to the identity.

\subsection{Group $\kerss$}

In this paragraph we prove that group $\kerss$ is isomorphic to
the group $\cH_{28}$ where
$$\cH_m\defin
\{
\begin{pmatrix}
1 & a & l\\
0 & 1 & b\\
0 & 0 & 1
\end{pmatrix}|a,b\in\zint~and~l\in\zint_m\}.$$
Group $\cH_m$ also can be described as the factor group of the
group $\cH$ (upper unitriangular $3\times 3$ matrices with integer
coefficients) by the cyclic subgroup generated by the matrix with
$a=b=0$ and $l=m$.

Idea of the proof is to compare presentations of two groups as we
just did above. First we give a presentation of $\kerss$. We know
already that $H^3(S^3\times S^3,\zint)$ is generated by
diffeomorphisms $A$ and $B$ (Lemma 2.). We also know from exact
sequence (\ref{secondseq}) that $\kerss$ is a central extension
of the group $\zint\oplus\zint$ by $\Theta_7\cong\faczin$ (see
\cite{KerMil}). If we denote by $\wH$ the generator of $\Theta_7$
(by the generator we mean a homotopy 7-sphere which bounds a
parallelizable manifold of signature 8) then we can choose $A,B$
and $\Sigma\defin\iota(\wH)$ as generators of $\kerss$. The
defining relations clearly will be: $A\Sigma=\Sigma
A,~B\Sigma=\Sigma B,~ABA^{-1}B^{-1}=\Sigma^k$, $\Sigma^{28}=1$ for
some $k\in\faczin$. So the goal is to figure out what this $k$ is.

We first define a map $\vars:~\kerss\lra\faczin$ as follows. Take
a representative $f\in{\it Diff(S^3\times S^3)}$ of a class
$[f]\in\kerss$.
\begin{dfn}
$$\vars([f])\defin D^4\times S^3~\bigcup_f~S^3\times D^4$$
Where $\bigcup_f$ means identification of a point $(x,y)\in
\partial(D^4\times S^3)$ with the point
$f(x,y)\in\partial(S^3\times D^4)$.
\end{dfn}

Denote by $\Sigma_f$ the manifold obtained from $\vars([f])$ by
smoothing the corners. It is clear that $\Sigma_f$ depends only
on the isotopy class $[f]$ and not on a specific representative
of this class. Note that $\Sigma_f$ is a homotopy sphere and
$\vars$ is a well defined map from $\kerss$ to the
group $\Theta_7\cong\zint_{28}$.\\

\noindent {\underline{Remark:}} Map $\vars$ is analog of the
Birman-Craggs homomorphism from the Torelli group of $S^2_g$
(2-dimensional surface of genus $g$) to $\zint_2$. See
\cite{Johnson} for the details.

\begin{lem}
The composition $\vars\circ\iota$ is the identity map of the group
$\Theta_7$.
\end{lem}
\begin{proof}
Take a sphere $\widetilde{\Sigma}_{\phi}\in\Theta_7$. We will
denote by $\phi$ the diffeomorphism of $S^3\times S^3$ which is
the identity outside an embedded disk $D^6\subset S^3\times S^3$
and corresponds to the element $\iota(\widetilde{\Sigma}_{\phi})$
of the mapping class group.

To show that $\vars\circ\iota=Id$ it is enough to show that
${\wS}_{\phi}$ is diffeomorphic to
$\Sigma_{\phi}=\vars\circ\iota(\wS_{\phi})$. We construct an
h-cobordism between these two manifolds. Take the handlebody
$D^4\times S^3$ and remove from it an interior disk $D^7$. The
resulting manifold is denoted by $\widetilde{D^4\times S^3}$.
Boundary components of $\widetilde{D^4\times S^3}$ are $S^6$ and
$S^3\times S^3$. Take disks $D^6_+$ in these two components and
connect them by a tube $D^6_+\times I$ embedded into
$\widetilde{D^4\times S^3}$. Next extend $\phi$ in an obvious way
(by the identity outside the tube) to a diffeomorphism $\Phi$ of
$\widetilde{D^4\times S^3}$. Consider now two manifolds:
$D^5\times S^3$ and $D^8$. Present the boundary of $D^5\times
S^3$ as the union:
$$
\partial(D^5\times S^3)=S^4\times S^3= D^4\times S^3\bigcup_
{S^3\times S^3}\widetilde{D^4\times S^3}\bigcup_{S^6}D^7
$$ and the boundary of $D^8$ as the union:
$$
\partial(D^8)=S^7=S^3\times D^4\bigcup_{S^3\times S^3}\widetilde{
D^4\times S^3}\bigcup_{S^6}D^7
$$
Using diffeomorphism $\Phi$ we can glue $D^5\times S^3$ and $D^8$
together along the common submanifold $\widetilde{ D^4\times
S^3}$ to obtain a cobordism (after smoothing the corners) $W^8$
between ${\widetilde{\Sigma}}_{\phi}$ and $\Sigma_{\phi}$. It is
clear that $W^8$ is simply connected. Using Mayer-Vietoris exact
sequence of the union $D^4\times S^3={\widetilde{D^4\times
S^3}}\cup D^7$ we see that
$$H_*(\widetilde{D^4\times S^3})\cong\left\{
\begin{array}{cc}
\zint & if~*=0,3,6\\
0 & otherwise
\end{array}\right.$$
In a similar way we can get homology groups of $W^8=D^5\times
S^3\cup D^8$:
$$
H_*(W^8)\cong\left\{
\begin{array}{cc}
\zint & if~*=0~or~7\\
0 & otherwise
\end{array}\right.~~~~and~~~~H_*(W^8,\Sigma_{\phi})\cong 0$$
Thus, by the h-cobordism theorem \cite{Milnor1} two homotopy
spheres ${\widetilde{\Sigma}}_{\phi}$ and $\Sigma_{\phi}$ are
diffeomorphic.
\end{proof}
\begin{thm}
The generators $A,B$ and $\Sigma$ of $\kerss$ satisfy the
relation: $ABA^{-1}B^{-1}=\Sigma$
\end{thm}
\begin{proof}
For the proof we construct a spin manifold $W^8$ bounded by the
sphere $\Sigma_{\comAB}$ and compute the $\mu$-invariant
$\mu(\Sigma_{\comAB})$ defined by Eells and Kuiper \cite{Eells}.

First we extend diffeomorphisms $A$ and $B$ to diffeomorphisms of
the handlebodies $D^4\times S^3$ and $S^3\times D^4$
respectively. It can be done since in the definition (recall
formulas (\ref{defAB})) we can assume that $s\in D^4,~t\in S^3$
for diffeomorphism $A$ and $s\in S^3,~t\in D^4$ for
diffeomorphism $B$. These extensions we also denote by $A$ and
$B$ respectively. Next we present $\Sigma_{\comAB}$ as the union
of five manifolds: $D^4\times S^3\bigcup_{\comAB}S^3\times
D^4~~\simeq$~
$$D^4\times S^3\bigcup_A S^3\times S^3\times
I\bigcup_B S^3\times S^3\times I\bigcup_{A^{-1}}S^3\times
S^3\times I\bigcup_{B^{-1}}S^3\times D^4$$ where $A,B,A^{-1}$ and
$B^{-1}$ belong to ${\it Diff(S^3\times S^3)}$. Consider
manifolds $D^8$, $D^5\times S^3$ and $S^3\times D^5$ with
boundaries presented as the unions:
$$\begin{array}{l}
\partial(D^5\times S^3)~=~S^4\times S^3~=~D^4\times
S^3\bigcup D^4\times S^3\\
\partial(S^3\times D^5)~=~S^3\times S^4~=~S^3\times D^4\bigcup
S^3\times D^4\\
\partial(D^8)~=~S^7~=~D^4\times S^3\bigcup S^3\times S^3\times
I\bigcup S^3\times D^4.
\end{array}$$
Using extension diffeomorphisms $A$ and $B$ defined above, we now
construct a manifold $W^8$ which will be used to compute
$\mu(\Sigma_{\comAB})$.
\begin{dfn}
Define $W^8$ to be the manifold obtained from the union
$$D^5\times S^3~\cup_A~D^8~\cup_B~
D^8~\cup_{A^{-1}}~D^8~\cup_{B^{-1}}~S^3\times D^5$$ by smoothing
the corners.
\end{dfn}
\begin{claim}
$D^5\times S^3\cup_A D^8~\simeq~D^8\cup_{B^{-1}}S^3\times
D^5~\simeq~D^8$
\end{claim}
\begin{proof}[Proof of the claim:] Evidently,
the union $D^5\times S^3\bigcup_{D^4\times S^3} D^8$ is a simply
connected manifold with simply connected boundary $D^4\times
S^3\bigcup_{S^3\times S^3}S^3\times D^4$. Using exact sequence of
Mayer-Vietoris it is easy to see that homology groups of
$D^5\times S^3\bigcup_{D^4\times S^3} D^8$ are trivial in all
dimensions $>0$. Hence by the {\it characterizations of the smooth
n-disk $D^n,~n\geq6$} (see \cite{Milnor1}) this union is
diffeomorphic to the disk $D^8$. Same proof works in the second
case.
\end{proof}
\noindent Thus we can write $W^8~=~D^8\cup_B D^8\cup_{A^{-1}}
D^8$. Now note that $\partial(W^8)=M(f_B,f_A)$ where (using
notations of Milnor \cite{Milnor2}) by $M(f_B,f_A)$ we denote the
boundary of the following union of three 8-disks:
$$(D^4\times D^4)_1\bigcup_{S^3\times D^4}(D^4\times
D^4)_2\bigcup_{D^4\times S^3}(D^4\times D^4)_3$$ The gluing maps
are (cf. \cite{Milnor2}, \S1):
$(x_1,y_1)\stackrel{f_B}{\lra}(x_2,y_2)\stackrel{f^{-1}_A}{\lra}
(x_3,y_3)$ where
$$y_3=y_2=f_B(x_1)\circ y_1,~~~~x_3=f_A(y_3)^{-1}\circ x_2=
f_A(y_3)^{-1}\circ x_1$$ and $f_B=f_A=f:~S^3\lra
SO(3)\stackrel{i_3}\hra SO(4)$, defined by the formula:
$f(x)\circ y=xyx^{-1}$ for $x\in S^3$ and $y\in SO(3)$. In
particular, we see that our homotopy sphere $\Sigma_{\comAB}$ is
diffeomorphic to the manifold $M(f_B,f_A)$. Using Mayer-Vietoris
exact sequence for the manifold $W^8\simeq
D^8\cup_{B}D^8\cup_{A^{-1}}D^8$ we see that
$H^*(W^8,\zint)\simeq0,~for~*=1,2~or~3$, and we can apply results
of Eells and Kuiper (\cite{Eells}, \S10) to the manifold $W^8$ to
compute the $\mu$-invariant of $M(f_B,f_A)\simeq\Sigma_{\comAB}$.
It is shown (see \cite{Eells}, page 109) that
$$
\mu(M(f_B,f_A))=\frac{B^2_1}{8(2!)^2}\left(1+\frac{2}{2^3-1}
\right)\left(\pm2p_1(f_B)p_1(f_A)\right)=\pm\frac{p_1(f_B)p_1(f_A)}{448}
$$
where $B_1=1/6$ is the first Bernoulli number and $p_1(f_B)$,
$p_1(f_A)$ are Pontrjagin numbers of the stable vector bundles
over $S^4$ determined by the compositions (cf. \cite{Milnor2},
\S3):
$$S^3\stackrel{f_B}\lra SO(4) \stackrel{i_4}\hookrightarrow
SO(5),~~~ S^3\stackrel{f_A}\lra SO(4)
\stackrel{i_4}\hookrightarrow SO(5)$$ We show that
$p_1(f_B)=p_1(f_A)=\pm4$. It is well known how Pontrjagin numbers
depend on the characteristic maps of the stable vector bundles
over the spheres : $p_s[S^{4s}]=\pm a_s\cdot\lambda\cdot(2s-1)!$.
Here $\lambda$ is the integer, corresponding to the
characteristic map: $S^{4s-1}\lra SO$ and $a_s=1~or~2$ if $s$ is
even or odd respectively. In our case $p_1(f_B)=\pm
2\cdot[i_4\circ f_B]$ and to find the integer $[i_4\circ f_B]$ we
need to find the homotopy class of the composition
$$S^3\lra SO(3)\stackrel{i_3}\hra SO(4)\stackrel{i_4}\hra SO(5).$$

It is a standard fact (which can be deduced from Theorem IV.1.12
of \cite{Browder2}) that inclusion $i=i_4\circ i_3$ induces map
$i_*: \pi_3(SO(3))\lra\pi_3(SO(5))$ which is multiplication by 2.
Hence $[i_4\circ f_B]=2$ and similarly $[i_4\circ f_A]=2$.
Therefore $\mu(\Sigma_{\comAB})=\mu(M(f_B,f_A))=
\pm\frac{16}{448}=\pm\frac{1}{28}$. For the generator $\wH$ of
$\zint_{28}$ we have (see \S4 of \cite{Eells}) $\mu(\wH)\equiv
\frac{-1}{2^57}\cdot8\pmod{1}$ or $\mu(\wH)=\frac{-1}{28}$. From
the theorem of Eells and Kuiper (\cite{Eells}, p.103) we see
(changing orientation of $\wH$ if necessary) that
$\Sigma_{\comAB}\simeq\wH$. Since $\comAB=\Sigma^k=\iota(\wH^k)$
we can apply Lemma 3. to get
$\wH^k=\vars\circ\iota(\wH^k)=\vars(\comAB)=\wH$. It shows that
$k\equiv 1\pmod{28}$ and finishes the proof.
\end{proof}
\noindent As a corollary we get a presentation of the group
$\kerss$:
\begin{equation}
\label{kernel} <A,B,\Sigma~|~\Sigma^{28}=1,~A\Sigma=\Sigma A,~
B\Sigma=\Sigma B,~\comAB=\Sigma>
\end{equation}
Consider now matrices:
$$A'= \begin{pmatrix}
1 & 1 & 0\\
0 & 1 & 0\\
0 & 0 & 1
\end{pmatrix};~~
B'= \begin{pmatrix}
1 & 0 & 0\\
0 & 1 & 1\\
0 & 0 & 1
\end{pmatrix};~~
\Sigma'= \begin{pmatrix}
1 & 0 & 1\\
0 & 1 & 0\\
0 & 0 & 1
\end{pmatrix}$$
It is easy to verify that these matrices $A',~B',~\Sigma'$
generate the group $\cH$ and satisfy the defining relations:
$<A',B',\Sigma'~|~A'\Sigma'=\Sigma'A',~
B'\Sigma'=\Sigma'B',~A'B'{A'}^{-1}{B'}^{-1}=\Sigma'>.$ Hence
groups $\kerss$ and $\cH_{28}$ are isomorphic.

\subsection{Generators and Relations}

We have shown in $\S 2.2$ that the factor group $\mcgss/\Theta_7$
admits the following presentation:
$<Y,U,A,B~|~YUY=UYU,~(YUY)^4=1,~AB=BA,~AY=YAB,~AU=UA,~BY=YB,~BU=UBA^{-1}>$.
Furthermore, we have shown explicitly (recall the proof of Theorem
1.) that $YUY=UYU,~AY=YAB,~AU=UA,~BY=YB,~BU=UBA^{-1}$ as
diffeomorphisms. Hence if we consider the group $\mcgss$ as a
central extension of $\mcgss/\Theta_7$ by the group $\Theta_7$, we
can choose $Y,U,A,B,\Sigma$ to be the generators. The defining
relations will be: $YUY=UYU,~ AY=YAB,~AU=UA,~BY=YB,~BU=UBA^{-1},~
~\Sigma^{28}=Id, ~\Sigma\leftrightarrows Y,U,A,B,$ and ~
$(YUY)^4=\Sigma^m,~ABA^{-1}B^{-1}=\Sigma^n$~ for some
$m,n\in\zint$. The symbol $\leftrightarrows$ means that element
on the left commutes with any element on the right. By Theorem 2.
we have $n=1$.
\begin{claim}
The following identities are valid in the group $\mcgss$:\\
$(B^{-1}YUY)^4=Id,~~YUYB^{-1}=A^{-1}\Sigma YUY,~~YUYA^{-1}=BYUY$.
\end{claim}
\begin{proof}
We know that $(B^{-1}YUY)^4=Id$ (Claim 1.). Using identities:
$AB=$ ~$BA\Sigma,~AY=YAB,~UB=BAU$ we find that
$YUYA^{-1}=YUA^{-1}B\Sigma Y=YA^{-1}UB\Sigma Y =YA^{-1}BA\Sigma
UY=YBUY=BYUY$. Similarly $YUYB^{-1}=
YUB^{-1}Y=YA^{-1}B^{-1}UY=A^{-1}B\Sigma B^{-1}YUY=A^{-1}\Sigma
YUY$.
\end{proof}

\noindent Now it is easy to see that $m=-1$. Indeed, \\
$Id=(B^{-1}~YUY)^2B^{-1}~YUYB^{-1}~YUY=
(B^{-1}~YUY)^2B^{-1}~ A^{-1}\Sigma (YUY)^2=$\\
$B^{-1}~ YUYB^{-1}~ YUY A^{-1}~B^{-1}(YUY)^2= B^{-1} YUYB^{-1}
BYUY B^{-1}(YUY)^2=$\\ $B^{-1} YUY A^{-1}\Sigma (YUY)^3=
B^{-1}B\Sigma (YUY)^4=\Sigma(YUY)^4$~~~ Hence
$(YUY)^4=\Sigma^{-1}$.

Let us now collect all the information we obtained so far and
state the main theorem of this paper.

\begin{thm}
The mapping class group of $S^3\times S^3$ admits the
following presentation:
$$\left\langle
\begin{array}{cc}
\Biggl.
\begin{array}{cc}
Y& U\\
A& B\\
\Sigma & ~
\end{array}
\Biggr| &
\begin{array}{ccc}
YUY=UYU,& (YUY)^4=\Sigma^{-1},& \Sigma\leftrightarrows Y,U,A,B \\
BU=UBA^{-1},& AY=YAB,& AB=BA\Sigma, \\
BY=YB,& AU=UA & ~
\end{array}
\end{array}
\right\rangle
$$
with
$Y=(y,0),~U=(u,0),~A=(a,0),~B=(b,0),~\Sigma=(id,1)\in\GJ\times\zint_{28}$.
\end{thm}

\noindent \underline{Remark:} It is well known that the mapping
class group of the 2-torus is generated by two Dehn twists.
Wajnryb \cite{Waj} has shown that for an orientable surface
$S^2_g$ of genus $g\geq2$ the mapping class group can be
generated by two elements that are not Dehn twists in general. It
follows from Theorem 3. and the work of Choie (see \cite{Choie},
Theorem 2.1.) that $\mcgss$ can also be generated by two elements.

\section{Cohomology of the Jacobi Group $\Gamma^J$}

It is usually difficult to obtain information about a group
having just generators and defining relations. The aim of this
part is to give an alternative description of the mapping class
group $\mcgss$ using the cohomology theory of groups. Since
$\mcgss$ is a central extension of the Jacobi group $\GJ$ by
$\zint_{28}$ it is natural to ask what element of the group
$H^2(\Gamma^J,\zint_{28})$ corresponds to this extension. First
we classify all central extensions of $\G1 $ by $\zint$ and
determine a 2-cocycle that generates $H^2(\G1 )$. In the second
section we show that $\GJ$ is isomorphic to an amalgamated
product. Using a Mayer-Vietoris exact sequence of this
amalgamated product we compute the cohomology groups of
$\Gamma^J$. Finally, we specify an element of
$H^2(\Gamma^J,\zint_{28})$ that corresponds to the mapping class
group $\mcgss$.

\subsection{Central extensions of $\G1 $ and $\ZZ$}

There are several ways to find cohomology groups of $\G1 $ with
trivial $\zint$-coefficients. All of these groups are of course
well known. The reason why we make some calculations here is that
in the next section these calculations will be generalized to
find the group $H^2(\GJ)$ and its generators.

Consider the following exact sequence: $ 0\lra \zint \lra E \lra
G \lra 0$. Group $E$ is called an extension of $G$ by
$\zint$. If the normal subgroup $\zint$ lies in the center of
$E$, this extension is called {\it central}. The equivalence classes
of such extensions are in 1-1 correspondence with the elements of the
second cohomology group $H^2(G,\zint)$. We will usually denote this
group by $H^2(G)$ forgetting the coefficients (only if $\zint$ is the
trivial $G$-module).

Recall that $\G1 \cong\zint_4*_{\zint_2}\zint_6$ and we can
consider a Mayer-Vietoris exact sequence:
\begin{equation}
\label{MayerGr} \ra H^n(\G1 )\ra H^n(\zint_4)\oplus
H^n(\zint_6)\ra H^n(\zint_2)\ra H^{n+1}(\G1 )\ra
\end{equation}

Cohomology groups of $\zint_m$ are known:
$H^{2k}(\zint_m)\cong\zint_m,~
H^{2k-1}(\zint_m)\cong0~for~k\geq1$. Hence $H^1(\G1 )\cong0$. We
also have the following fragments:
$$ 0\ra H^{2n}(\G1 )\ra H^{2n}(\zint_4)\oplus H^{2n}(\zint_6)
\stackrel{j^*}{\ra} H^{2n}(\zint_2)\ra H^{2n+1}(\G1 )\ra 0$$ with
$j^*=i^*_4+i^*_6$ and $i_k:\zint_2\lra\zint_k$ - multiplication
by $k/2$. If we denote by $t,~z$ and $z'$ generators of groups
$H^2(\zint_2),~H^2(\zint_4)$ and $H^2(\zint_6)$ respectively,
then it is easy to see that $i^*_4(z)=i^*_6(z')=t$ i.e.
$j^*(n,m)=n+m$. Hence $Im(j^*)=H^2(\zint_2)$, $Ker(j^*)$ is
generated by the element $(1,5)\in \zint_4\oplus\zint_6$ and
isomorphic to $\zint_{12}$. Thus $H^3(\G1 )\cong0$ and $H^2(\G1
)\cong\zint_{12}$. It follows from the properties of the {\it
norm map} (\cite{Brown}, ch.III) that $H^{2k-1} (\G1 )\cong0$ and
$H^{2k}(\G1 )\cong\zint_{12}$ for any $k\geq1$.

Now we write down an explicit function $f: \G1 \times \G1
\lra\zint$ that generates the group $H^2(\G1 )$. This function
will be defined in terms of a generator of the group
$H^2(\zint_{12})$ which we describe first. Consider the group
$\zint_m$ as the multiplicative group on elements
$\{z^1,\dots,~z^p,\dots,~z^{m-1},id\}$ and define a function $f_m:
\zint_m\times\zint_m\lra\zint$ by the formula:
\begin{dfn}
\begin{equation}
\label{generatorZm} f_m(z^p,z^q)\defin\left\{
\begin{array}{ccc}
1 & if & \bar p+\bar q\geq m\\
0 & if & \bar p+\bar q < m
\end{array}
\right.
\end{equation}
where $\bar p,\bar q\in \{0,1,2,\dots, m-1\}$ and $\bar p\equiv p
\pmod{m},~\bar q\equiv q \pmod{m}$
\end{dfn}
\begin{lem}
Function $f_m$ is a generator of the group
$H^2(\zint_m)\cong\zint_m$.
\end{lem}
\begin{proof}
Equality $f_m(z^p,z^q)+f_m(z^{p+q},z^r)=
f_m(z^q,z^r)+f_m(z^p,z^{q+r})$ shows that $f_m$ is a 2-cocycle.
Verification of this equality is straightforward and left as an
exercise. Let $<c>$ and $<z~|~z^m=id>$ be presentations of groups
$\zint$ and $\zint_m$ respectively. Then a central extension of
$\zint_m$ by $\zint$ will have a presentation:
$<Z,C~|~ZC=CZ,~Z^m=C^k>$ with $Z=(z,0),~C=(id,1)
\in\zint_m\times\zint$ and some $k\in\zint$. It follows from
formula (\ref{generatorZm}) that $k=1$ for the function $f_m$,
and that the cocycle $t\cdot f_m$ defines the extension with
$k=t$. We denote this extension by $E_t$. If $m\nmid s-t$ then
cocycles $s\cdot f_m$ and $t\cdot f_m$ define non equivalent
extensions. Indeed, suppose that $E_s$ and $E_t$ are equivalent.
Then there exists an isomorphism $\iota: E_s\lra E_t$ that makes
the following diagram commute:
$$\begin{CD}
0 @>>> \zint @>i_s>> E_s @>\rho_s>> \zint_m @>>> 0\\
 @.    @VV=V       @VV\iota V      @VV=V   @.\\
0 @>>> \zint @>i_t>> E_t @>\rho_t>> \zint_m @>>> 0
\end{CD}$$
Using the above presentation of groups $E_s$ and $E_t$ one can see
from this diagram that $s-t=lm$ for some $l\in\zint$. Hence $f_m$
is a generator of the group $H^2(\zint_m)$.
\end{proof}

Now consider the canonical projection $ab: \G1 \lra (\G1 )_{ab}
\cong\zint_{12}$. We define the function $f:\G1 \times \G1
\lra\zint$ by the formula
\begin{dfn}
$f(M,N)\defin f_{12}(M_{ab},N_{ab})$ where
$M_{ab}=ab(M)\in\zint_{12}$
\end{dfn}
\begin{lem}
Function $f(M,N)=f_{12}(M_{ab},N_{ab})$ is a generator of the
group $H^2(\G1 )$.
\end{lem}
\begin{proof}
Since $(M\cdot N)_{ab}=M_{ab}\cdot N_{ab}$ it follows from the
previous lemma that $f$ is a 2-cocycle. Let $<c>$ and $<y,u~|~
yuy=uyu,~(yuy)^4=id>$ be presentations of groups $\zint$ and $\G1
$. Then a central extension of $\G1 $ has a presentation:
$$<Y,U,C~|~YUY=UYUC^a,~(YUY)^4=C^b,~CY=YC,~CU=UC>$$
with $Y=(y,0),~ U=(u,0),~C=(id,1)$ and some $a,b\in\zint$. Using
the group law: $(g,k)\cdot(h,l)=(gh,k+l+f(g,h))$ one can easily
find that function $f$ defines the extension with $a=0$ and $b=1$.
Cocycle $t\cdot f$ defines the extension with $a=0$ and $b=t$. As
in the proof of Lemma 4. one can show that if $12\nmid s-t$ then
cocycles $s\cdot f$ and $t\cdot f$ define non equivalent
extensions. Therefore function $f$ is indeed a generator of
$H^2(\G1 )$.
\end{proof}

Since the 2-torus is an {\it Eilenberg-MacLane complex} of type
$(\zint^2,1)$ we have $H^2(\zint^2)\cong\zint$. We can write an
explicit function $g:\zint^2 \times\zint^2\lra \zint$ that
generates this cohomology group. For instance, take two vectors
$v_1=(a,b),~v_2=(s,t)\in\zint^2$ and consider
$g:~\zint^2\times\zint^2\lra\zint$ defined by the formula:
$g(v_1,v_2)\defin at$. One can check that $g$ is a 2-cocycle which
defines the following extension:
\begin{equation}
\label{zz} <A,B,C~|~AC=CA,~BC=CB,~ABA^{-1}B^{-1}=C>
\end{equation}
with $A=((1,0);0),~B=((0,1);0),~C=((0,0);1)\in
\zint^2\times\zint$. It can be shown (as we did in the proof of
Lemma 4.) that any function which defines extension (\ref{zz}) is
a generator of the group $H^2(\zint^2)$.

\noindent {\underline{Remark:}} Consider another function
$\varphi:~\zint^2\times\zint^2\lra\zint$ defined as follows. For
$v=(\lambda,\mu)$, and
$w=(\lambda',\mu'),~~\varphi(v,w)\defin\begin{vmatrix} \lambda
&\mu \\\lambda' &\mu'\end{vmatrix}=\lambda\mu'-\lambda'\mu$. It
can be easily verified that $\varphi$ is a cocycle cohomologous
to $2g$. Note that $\varphi$ defines the extension which is
isomorphic to the {\it Heisenberg group}~~$H(\zint)$. This later
one is~ $\zint\oplus\zint\oplus\zint$~ as a set and with
multiplication:
$$(a,b,c)\cdot(a',b',c')\defin (a+a',b+b',c+c'+ a\cdot b'- b\cdot
a')$$ Here one can consider so called (real) Jacobi group
$G^J(\zint)$ (cf. \cite{Berndt}) defined as the semidirect
product of $\G1 $ and $H(\zint)$: $G^J(\zint)=\G1
\ltimes H(\zint)$.

\subsection{Group $\Gamma^J$ as an amalgamated product}

In this paragraph we will compute the cohomology groups of $\GJ$.
First we present $\GJ$ as an amalgamated product and then use a
Mayer-Vietoris exact sequence of this product to find $H^2(\GJ)$.

Consider for $m=2,4,6$ the cyclic groups $\zint_m$ generated by
the matrices $\begin{pmatrix} -1&0\\0&-1\end{pmatrix}$,
$\begin{pmatrix} 0&-1\\1&0\end{pmatrix}$, $\begin{pmatrix}
0&1\\-1&1\end{pmatrix}$ respectively. These are subgroups of $\G1
$ and they act on the elements of $\ZZ$ in the same natural way
as $\G1 $ does. With respect to this action we define the
semidirect products $G_m\defin \zint_m\ltimes\zint^2$.

\begin{lem}
Group $\GJ$ is isomorphic to the amalgamated product
$G_4*_{G_2}G_6$.
\end{lem}
\begin{proof}
First we give presentations of groups $G_m$ with generators and
defining relations. Denote the generators of $\zint_2,~\zint_4$,
$\zint_6$ by $\alpha,~\beta$, $\gamma$ and elements $(id,(1,0))$,
$(id,(0,1))$ of the group $G_m$ by $A$ and $B$
respectively.\\
{\underline{Note:}} To avoid cumbersome notations we use the same
letters $A$ and $B$ to denote elements of different groups. We
hope it will not cause any confusion.\\
Using group law $(M,X)\cdot(M',X')= (MM',XM'+X)$ one can easily
show (cf. proof of Lemma 1. above) that $A\alpha=\alpha
A^{-1},~B\alpha=\alpha B^{-1},~B\beta= \beta A^{-1},~A\beta=\beta
B$ and $A\gamma=\gamma B,~B\gamma=A\gamma
A^{-1}$. Hence we get presentations:\\
$G_2\simeq <A,B,\alpha~|~AB=BA,~\alpha^2 =id,~B\alpha=\alpha
B^{-1},~A\alpha=\alpha A^{-1}>$\\ $\Rightarrow~ (G_2)_{ab}\simeq
\zint_2\oplus \zint_2\oplus\zint_2$,\\
$G_4\simeq <A,B,\beta~|~AB=BA,~\beta^4=id,~B\beta= \beta
A^{-1},~A\beta=\beta B>$\\ $\Rightarrow~ (G_4)_{ab}\simeq
\zint_4\oplus\zint_2$,\\
$G_6\simeq <A,B,\gamma~|~AB=BA,~\gamma^6=id,~B\gamma=A\gamma
A^{-1},~A\gamma=\gamma B>$\\ $\Rightarrow~ (G_6)_{ab}\simeq
\zint_6$.

Next define maps $\iota_4:~G_2\lra G_4$ and $\iota_6:~G_2\lra G_6$
by the formulas: $\iota_4(A)=A$, $\iota_4(B)=B$,
$\iota_4(\alpha)=\beta^2$ and $\iota_6(A)=A,~\iota_6(B)=B,~\iota_6
(\alpha)=\gamma^3$. Map $\iota_4$ induces commutative diagram
\begin{equation}
\label{G4}
\begin{CD}
0 @>>> \zint^2 @>>> G_2 @>>> \zint_2 @>>> 0\\
@.    @VV{\cong}V   @VV{\iota_4}V  @VV{\times2}V  @.\\
0 @>>> \zint^2 @>>> G_4 @>>> \zint_4 @>>> 0
\end{CD}
\end{equation}
from which it follows that $\iota_4$ is a monomorphism. Similarly
one proves that $\iota_6$ is a monomorphism. Since $\alpha$ is
identified with $\beta^2$ and $\gamma^3$ we obtain the following
presentation of $G_4*_{G_2}G_6$ with generators and defining
relations:
$$\left\langle
\begin{array}{cc}
\Biggl.
\begin{array}{cc}
\beta& \gamma \\
A& B
\end{array}
\Biggr| &
\begin{array}{cccc}
AB=BA, & A\beta=\beta B, & A\gamma=\gamma B, & B\gamma=A\gamma A^{-1}\\
B\beta=\beta A^{-1}, & \beta^2=\gamma^3, & \beta^4=id & ~
\end{array}
\end{array}
\right\rangle
$$

Consider now two elements: $U\defin \beta\gamma^{-1}$ and
$Y\defin\gamma^2\beta^{-1}$. Obviously, $YU=\gamma$ and
$UYU=\beta$. One can easily show that the above presentation of
the group $G_4*_{G_2}G_6$ is equivalent to one of the Jacobi
group $\Gamma^J$ obtained in Lemma 1. \end{proof}

To find the cohomology of $G_m$ we will use the
Lyndon-Hochschild-Serre (LHS) spectral sequence (see \cite{Hoch}
or \cite{Evens}, \S7.2) of the split extension that defines $G_m$:
$E_2^{p,q}=H^p(\zint_m,H^q(\zint^2,\zint))\Rightarrow H^{p+q}
(G_m,\zint)$. We need to know how $\zint_m$ acts on
$H^q(\zint^2)$. Suppose $f\in H^q(\zint^2)$, that is
$f:\zint^2\times\cdots\times\zint^2 \lra \zint$. Then we have
$\cM\circ f(\sigma_1,\dots,\sigma_q)= f(\cM^{-1}\sigma_1\cM,\dots,
\cM^{-1}\sigma_q\cM)$ (since $\zint$ is a trivial
$\zint_m$-module) where $\sigma_i=(\begin{pmatrix}
1&0\\0&1\end{pmatrix};(x_i,y_i))$ and $\cM=(\begin{pmatrix}
a&b\\c&d\end{pmatrix};(0,0))$ are the corresponding elements of
the group $G_m$ (cf. \cite{Hoch}, p.117). If we denote by $M$ the
matrix of the element $\cM$ and by $\sigma_i$ the vector of the
element $\sigma_i$, we find that $M\circ
f(\sigma_i,\dots,\sigma_q) = f(\sigma_i\cdot M,\dots,\sigma_q\cdot
M)$ where on the right we mean multiplication of the row vector
by a matrix. Groups $H^q(\zint^2)$ are nonzero only for $q=0,1,2$.

It is easy to verify that $\zint_m$ acts trivially on groups
$H^0(\zint^2)\cong\zint$ and $H^2(\zint^2)\cong\zint$. If we
denote by $(n,k)$ an element of $H^1(\zint^2)\cong\ZZ$, it can be
shown that $M\circ(n,k)=(n,k)\cdot M^T$.

Since groups $H^*(\zint_m)$ are well known we only need to
compute $H^*(\zint_m,\zint^2)$ with respect to the action
described above. $H^*(\zint_m,\zint^2)$ can be found using the
{\it norm map} $\bar N:~M_G\lra M^G$ for group $G=\zint_m$ and
the module $M=\zint^2$, since $H^{2k}(G,M)\cong coker \bar N$ and
$H^{2k-1}(G,M)\cong ker\bar N$ for any $k\geq1$ (see \cite{Brown},
ch.III).

\noindent 1) For $\zint_2$ we have $M_G\cong
\zint_2\oplus\zint_2$ with the generators $(1,0),(0,1)$. It is obvious
that $M^G\cong 0$ and we get:
$$H^{2k}(\zint_2,\zint^2)\cong0,~~~~~H^{2k+1}(\zint_2,\zint^2)
\cong\zint_2\oplus\zint_2~~~~~ for ~any~k\geq0$$

\noindent 2) For $\zint_4$ we have $M_G\cong\zint_2$ with the
generator $(1,0)$. Here $M^G\cong 0$ also, and we get:
$H^{2k}(\zint_4,\zint^2)\cong0,~~~~~H^{2k+1}(\zint_4,\zint^2)
\cong\zint_2~~~~~ for ~any~k\geq0$.

\noindent 3) In the case of $\zint_6$, $M_G\cong M^G\cong 0$ and
$H^k(\zint_6,\zint^2)\cong0~~~~~ for ~any~ k\geq 0$.

\begin{lem} For any $k\geq 0$ and $m\in\{2,4,6\}$,
$H^{2k+1}(G_m)\cong 0$. Furthermore,\\
$H^2(G_2)\cong~H^0(\zint_2)\oplus H^1(\zint_2,\ZZ)\oplus
H^2(\zint_2)~\cong~\zint\oplus(\zint_2\oplus\zint_2)\oplus\zint_2$,\\
$H^2(G_4)\cong~H^0(\zint_4)\oplus H^1(\zint_4,\ZZ)\oplus
H^2(\zint_4)~\cong~\zint\oplus\zint_2\oplus\zint_4$,\\
$H^2(G_6)\cong~H^0(\zint_6)\oplus H^2(\zint_6)~\cong~
\zint\oplus\zint_6$.
\end{lem}
\begin{proof}
We prove this lemma only for $m=2$. In the other cases the proof
is similar. Consider the first quadrant of the LHS spectral
sequence $E_2^{p,q}=H^p(\zint_2,H^q(\zint^2))$, for
$p,q\geq0$.\\
$$
\underrightarrow{\left\uparrow
\begin{array}{ccccccc}
\vdots &\vdots &\vdots &\vdots &\vdots &\vdots & ~ \\
0 &0 &0 &0 &0 &0 &\ldots\\
\zint &0 &\zint_2 &0 &\zint_2 &0 &\ldots\\
0 &\zint_2^2 &0 &\zint_2^2 &0 &\zint_2^2 &\ldots\\
\zint &0 &\zint_2 &0 &\zint_2 &0 &\ldots
\end{array}\right.}
$$\\
From this $E_2$-term we see that $E^{*,*}_2 \cong
E^{*,*}_3\cong...\cong E^{*,*}_{\infty}$. Hence for any $k\geq
0$, $H^{2k+1}(G_2)\cong 0$. Since
$H_1(G_2)\cong(G_2)_{ab}\cong\zint_2\oplus\zint_2\oplus\zint_2$,
by the Universal Coefficient Theorem $H^2(G_2)\cong {\it
Hom}(H_2(G_2),\zint)\oplus{\it Ext} (H_1(G_2),\zint)\cong{\it
Hom}(H_2(G_2),\zint)\oplus\zint_2^3$, and therefore
$H^2(G_2)\cong H^0(\zint_2)\oplus H^1(\zint_2,\zint^2) \oplus
H^2(\zint_2)~\cong~\zint\oplus
(\zint_2\oplus\zint_2)\oplus\zint_2$.
\end{proof}
\begin{thm}
$\GJ$ has the following homology and cohomology groups:\\
$H_1(\GJ)\cong\zint_{12},~~~H_2(\GJ)\cong\zint\oplus\zint_2,~~~
H^1(\GJ)\cong 0,~~~H^2(\GJ)\cong\zint\oplus\zint_{12}$,
$H^3(\GJ)\cong\zint_2$
\end{thm}
\begin{proof}
First note that it follows from the presentation of the group
$\GJ$ that $H_1(\GJ)\cong(\GJ)_{ab}\cong\zint_{12}$. To
compute the cohomology groups, we use a Mayer-Vietoris exact
sequence of the amalgamated product $G_4 *_{G_2} G_6$. By the
previous lemma we get the following fragment of this exact
sequence:
$$0 \lra H^2(\GJ)\lra H^2(G_4)\oplus
H^2(G_6)\stackrel{j}{\lra} H^2(G_2)\lra H^3(\GJ)\lra 0$$ where
$j=\iota_4^*+\iota_6^*$ and $\iota_m^*$ is induced by the
inclusion $\iota_m:G_2\hra G_m$ (cf. proof of Lemma 6.).
Commutative diagram (\ref{G4}) together with the {\it functorial
dependence} (cf. \cite{Evens}, \S7.2) induces the functorial map~
$H^p(\zint_4,H^q(\zint^2))\lra H^p(\zint_2,H^q(\zint^2))$. Hence
$H^2(G_4)\stackrel{\iota_4^*}{\lra}H^2(G_2)$ is the map induced by
multiplication: $\zint_2\stackrel{\times2}{\lra}\zint_4$. If we
denote by $x,y,z$ generators of groups $H^0(\zint_4)\cong\zint$,
$H^1(\zint_4,\zint^2)\cong\zint_2 ,~H^2(\zint_4)\cong\zint_4$
respectively then we see that $\iota_4^*(x,y,z)=(x',y',0,z')$
where $x',y',z'$ are generators of the groups $\zint\cong
H^0(\zint_2),~\zint_2\subset H^1(\zint_2,\zint^2),~\zint_2\cong
H^2(\zint_2)$. Similarly one can show that
$\iota_6^*(s,t)=(x',0,0,z')$ where by $s$ and $t$ we denoted the
generators of groups $H^0(\zint_6)$ and $H^2(\zint_6)$. Hence
$H^3(\GJ)\cong Coker(j)\cong\zint_2$ and $H^2(\GJ)\cong
Ker(j)\cong \zint\oplus\zint_{12}$. It is clear that
$H^1(\GJ)\cong0$. We can assume that
$H_2(\GJ)\cong\zint^{n_2}\oplus Tor_2$ and
$H_3(\GJ)\cong\zint^{n_3}\oplus Tor_3$ where by $Tor_k$ we denote
the elements of finite order. Then by the Universal Coefficient
Theorem $\zint\oplus\zint_{12}\cong H^2(\GJ)\cong {\it
Ext}(\zint_{12},\zint)\oplus{\it Hom}(\zint^{n_2}\oplus
Tor_2,\zint)\cong\zint_{12}\oplus\zint^{n_2}$, that is $n_2=1$.
Similarly we have $\zint_2\cong H^3(\GJ)\cong{\it Ext(\zint\oplus
Tor_2,\zint)\oplus Hom(\zint^{n_3}\oplus Tor_3,\zint)}\cong{\it
Ext(Tor_2,\zint)}\oplus\zint^{n_3}$. Therefore $n_3=0$ and
$Tor_2\cong\zint_2$ as required.
\end{proof}

\noindent {\underline{Remark:}} To find the other cohomology
groups of $\GJ$ one could use the LHS spectral sequence of the
defining extension of $\GJ$: $H^p(\G1 ,H^q(\zint^2)) \Rightarrow
H^{p+q}(\GJ)$. It can be seen that $H^*(\G1 ,H^2(\zint^2))\cong
H^*(\G1 ,H^0(\zint^2))\cong H^*(\zint_{12})$; $H^{2k}(\G1
,H^1(\zint^2))\cong \zint_2$ and $H^{2k-1}(\G1
,H^1(\zint^2))\cong 0$ for any $k\geq1$. The previous theorem
implies that $d_2^{0,2}$ of this spectral sequence is the zero
map. Using the cup product we deduce that $d_2^{2,2}$ is zero,
hence $d_2^{4,2}$ is zero and so on. Applying Proposition 7.3.2.
of \cite{Evens} and Theorem 4. of \cite{Hoch} we see that
$$
H^{2k+1}(\GJ)\cong\zint_2~~and~~H^{2k+2}(\GJ)\cong\zint_{12}
\oplus\zint_{12}~~~~~for~any~k\geq1$$

\begin{thm}
The split extension $0\lra\zint^2 \stackrel{i}{\lra}\GJ
\stackrel{\rho}{\lra}\G1 \lra 0$, induces the following split
exact sequence (compare with the five-term exact sequence
\cite{Evens}, \S 7.2):
\begin{equation}
\label{split}
0\lra H^2(\G1 )\stackrel{\rho^*}{\lra} H^2(\GJ)
\stackrel{i^*} {\lra} H^2(\zint^2) \lra 0
\end{equation}
\end{thm}
\begin{proof}
Choose the presentation of $\GJ$ found in Lemma 1. Then any
central extension $E$: $0\ra\zint\ra E\stackrel{\pi}{\ra}
\GJ\ra0$  has the following presentation:
$$
\left\langle
\begin{array}{cc}
\Biggl.
\begin{array}{cc}
Y& U\\
A& B\\
\Sigma & ~
\end{array}
\Biggr| &
\begin{array}{ccc}
YUY=UYU\Sigma^{k_0},& (YUY)^4=\Sigma^{k_1},& AB=BA\Sigma^{k_2},\\
BU=UBA^{-1}\Sigma^{k_3},& AY=YAB\Sigma^{k_4},& AU=UA\Sigma^{k_5},\\
BY=YB\Sigma^{k_6},& \Sigma\leftrightarrows Y,U,A,B & ~
\end{array}
\end{array}
\right\rangle
$$
with some $k_0,\dots,k_6\in\zint$ and
$Y=(y,0),~U=(u,0),~A=(a,0)$, $B=(b,0)$ and $\Sigma=(id,1)$.
Consider element ${\bf U}\defin U\Sigma^{k_0}$. Then from equality
$YUY=UYU\Sigma^{k_0}$ we get $Y{\bf U}Y={\bf U}Y{\bf U}$. Instead
of $(YUY)^4=\Sigma^{k_1}$ we get $(Y{\bf U}Y)^4=
\Sigma^{k_1+4k_0}$ and so on. Thus, by changing $U$ to
$U\Sigma^{k_0}$ we can make $k_0=0$. Similarly, changing $A$ to
$A\Sigma^{-k_3}$ and $B$ to $B\Sigma^{k_4}$ we can eliminate $k_3$
and $k_4$ and assume that $k_0=k_3=k_4=0$ in the above
presentation of $E$. Now one can easily obtain the following
equalities: $BYUY=YUYA^{-1}\Sigma^{2k_6}$,
$BUYU=UYUA^{-1}\Sigma^{k_6-k_5}$. Therefore $k_6=-k_5$.
Analogously, if we compare $AYUY$ with $AUYU$ we find that
$k_6=k_5$. Hence $k_5=k_6=0$. It means that we can choose
set-theoretic cross-sections of $\pi$ (i.e., functions
$s_{m,n}:\GJ\lra E$ so that $\pi\circ s_{m,n}=id$) in such a way
that any central extension of $\GJ$ by $\zint$ has a presentation:
$$
E_{m,n}\defin
\left\langle
\begin{array}{cc}
\Biggl.
\begin{array}{cc}
Y& U\\
A& B\\
\Sigma & ~
\end{array}
\Biggr| &
\begin{array}{ccc}
YUY=UYU,& (YUY)^4=\Sigma^n,& AB=BA\Sigma^m,\\
BU=UBA^{-1},& AY=YAB,& AU=UA,\\
BY=YB,& \Sigma\leftrightarrows Y,U,A,B & ~
\end{array}
\end{array}
\right\rangle
$$
and any element of $H^2(\GJ)$ is cohomologous to a cocycle,
defined by one of these $s_{m,n}$.

Choose a 2-cocycle $\omega_1$ that defines the extension
$E_{1,0}$. Evidently, $i^*(\omega_1)$ defines the extension
(\ref{zz}) of $\zint^2$ by $\zint$ and therefore generates
$H^2(\zint^2)$. Hence $i^*$ is onto and $\omega_1$ generates the
direct summand $\zint$ of $H^2(\GJ)$. As for the other generator
$\omega_2$ of $\zint_{12}\subset H^2(\GJ)$, we can choose
$\rho^*(f)$ (recall Def. 4 above). It can be verified that
$\omega_2=\rho^*(f)$ defines the extension $E_{0,1}$ and the
subgroup generated by $\omega_2$ is the kernel of $i^*$. This
proves the exactness. Proof of the asserted splitting is left to
the reader (cf. Proposition 7.3.2. of \cite{Evens}).
\end{proof}

\subsection{On groups $H^2(\GJ,\zint_{28})$ and $\mcgss$}

In this short paragraph we determine an element of
$H^2(\GJ,\zint_{28})$ that corresponds to the mapping class group
of $S^3\times S^3$.

Consider a short exact sequence $0\lra A\stackrel{\mu}{\lra} B
\stackrel{\nu}{\lra} C\lra 0$ of $G$-modules. Then we know from
homological algebra that there is a natural map $\delta:
H^n(G,C)\lra H^{n+1}(G,A)$ such that the sequence
$$
\dots\lra H^n(G,A)\stackrel{\mu_n}{\lra}
H^n(G,B)\stackrel{\nu_n}{\lra} H^n(G,C)\stackrel{\delta}{\lra}
H^{n+1}(G,A)\lra\dots
$$
is exact. Therefore the short exact sequence $0\lra
\zint\stackrel{\times 28}{\lra} \zint \lra \zint_{28}\lra 0$ of
trivial $\GJ$-modules gives the long exact sequence
$$
\lra H^2(\GJ)\stackrel{\mu_2}{\lra} H^2(\GJ)\stackrel{\nu_2}{\lra}
H^2(\GJ,\zint_{28})\stackrel{\delta}{\lra}
H^3(\GJ)\stackrel{\mu_3}{\lra} H^3(\GJ) \lra
$$
Map $\mu_n$ is multiplication by 28 and we have (from Theorem 4.)
$H^2(\GJ,\zint_{28})\cong Im(\nu_2)\oplus Ker(\mu_3)\cong
\zint_{28}\oplus\zint_4\oplus\zint_2$. It follows from the
presentation of $\mcgss$ given in Theorem 3. that this mapping
class group is the factor group of the central extension of $\GJ$
that corresponds to cocycle $\omega_1-\omega_2$ of $H^2(\GJ)$. It
is clear that $\nu_2(\omega_1-\omega_2)$ will be a cocycle that
corresponds to the group $\mcgss$. If we denote a generator of
the summand $\zint_2$ by $\to_3$, and generators of
$\zint_{28},~\zint_4$ by $\to_1\defin\nu_2(\omega_1)$ and
$\to_2\defin\nu_2(\omega_2)$ respectively, then these three
cocycles $\to_i,~i\in\{1,2,3\}$ generate group
$H^2(\GJ,\zint_{28})$ and $\to_1-\to_2$ will be a cocycle that
defines an extension isomorphic to $\mcgss$.

\noindent {\underline{Remark:}} One can deduce from the proof of
Theorem 5. that cocycle $\to_3$ defines the extension $E$ with
$k_0=k_1=k_2=k_3=k_4=0$ and $k_5=k_6=1$.

\end{document}